\documentclass[11pt]{amsart}
\usepackage{amsmath,amsfonts,amsthm,amssymb,nameref}
\usepackage[top=4cm, bottom=3cm, left=3.3cm, right=3.3cm]{geometry} 

\usepackage{tikz} 
\usetikzlibrary{calc,automata} 

\usepackage{color} 

\usepackage{enumerate}

\usepackage[pdftex]{hyperref}
\hypersetup{colorlinks,citecolor=black,filecolor=black,linkcolor=black,urlcolor=black}
\usepackage[all]{hypcap}
\usepackage{microtype} 

\newcommand{\N}	{\mathbb N}

\newcommand{\matP}{\mathbb P}

\newcommand{\Cay}	{\operatorname{Cay}}
\newcommand{\diam}	{\operatorname{diam}} 
\newcommand{\girth}	{\operatorname{girth}} 

\newtheorem{thm}{Theorem}

\newtheorem*{thm*} {Theorem} 
\newtheorem*{prop*}{Proposition}
\newtheorem*{lem*} {Lemma} 
\newtheorem*{cor*} {Corollary}
\theoremstyle{definition}

\newtheorem{remark}[thm]{Remark}
\newtheorem*{defi*}{Definition}
\newtheorem*{example*}{Example}
\newtheorem*{remark*}{Remark}
\newtheorem*{problem*}{Problem}
\newtheorem*{convention*}{Convention}

\begin{document}
\title[Gromov's monsters do not act on hyperbolic spaces]{Gromov's random monsters do not act non-elementarily on hyperbolic spaces}
\author{Dominik Gruber}
\address{Department of Mathematics, ETH Zurich, 8092 Zurich, Switzerland}
 \email{dominik.gruber@math.ethz.ch}

\author{Alessandro Sisto}
\address{Department of Mathematics, ETH Zurich, 8092 Zurich, Switzerland}
 \email{sisto@math.ethz.ch}
 
 \author{Romain Tessera}
\address{Laboratoire de Math\'ematiques d'Orsay, Univ. Paris-Sud, CNRS, Universit\'e Paris-Saclay, F-91405 Orsay, France}
\email{tessera@phare.normalesup.org}

 \begin{abstract}
We show that Gromov's monster groups arising from i.i.d.\ labelings of expander graphs do not admit non-elementary actions on geodesic hyperbolic spaces.
The proof relies on comparing properties of random walks on randomly labeled graphs and on groups acting non-elementarily on hyperbolic spaces.
\end{abstract}

\maketitle

\section*{Introduction}

In ``Spaces and questions" \cite{Gromov_questions},  Gromov introduced a  random model of finitely generated infinitely presented groups. The main motivation was to prove the existence of finitely generated groups with extreme geometrical properties. Let us briefly outline the construction. Given three integers $d\geq 3$, $k\geq 2$ and $j\geq 1$ and an infinite sequence of finite $d$-regular graphs $(\Omega_n)$, we label independently at random each edge of every $\Omega_n$ with a word of length $j$ in the free group on $k$ generators $F_k$. We then define a group $G$ as the quotient of $F_k$ by the normal closure of the set of words corresponding to closed loops of the $\Omega_n$. The main interest of this construction is that for appropriate parameters $d,k,j$, if one choses the $\Omega_n$ to be a suitable sequence of expander graphs, the group $G$ does not coarsely embed into Hilbert space, because it contains an \emph{almost quasi-isometrically} embedded copy of the sequence $(\Omega_n)$ in its Cayley graph \cite{Gromov_random,Arzhantseva-Delzant}.  Moreover, it provides a counterexample to the Baum-Connes conjecture with coefficients \cite{HLS}.

The main result of this note is that, under very general conditions on the graphs $\Omega_n$ (which, in particular, are satisfied by the expanders considered by Gromov), the group $G$ almost surely does not act non-elementarily\footnote{By ``non-elementary action" we mean possessing two hyperbolic elements with disjoint pairs of fixed points on the boundary.} on any geodesic Gromov hyperbolic space.

Our motivation for proving this fact comes from another type of Gromov's monsters whose existence has been proven in \cite{Osa-label}, where a counting argument is used to show the existence of labelings on the $\Omega_n$ satisfying a certain ``small cancellation" property (which is not satisfied by the i.i.d.\ labeling described above). They satisfy the graphical small cancellation condition developed in \cite{Gromov_random,Ollivier}, see also \cite{Gru-TAMS} (as opposed to the \emph{geometric} small cancellation condition satisfied by the i.i.d.\ labelings \cite{Gromov_random,Arzhantseva-Delzant}). This graphical small cancellation condition ensures that the resulting Cayley graph contains an {\it isometric} copy of each $\Omega_n$.
All (infinitely presented) graphical small cancellation groups are {\it acylindrically hyperbolic} \cite{GS-smallcanc}, and therefore, they admit non-elementary actions on geodesic hyperbolic spaces (even on quasi-trees \cite{Balas}). Thus, our result provides a strong distinction between the two types of Gromov's monster groups.

We explain the reason for the sharp difference between the groups coming from the i.i.d.\ labeling and those coming from the graphical small cancellation labeling.

Our argument combines crucial estimates on random walks on $G$ of Naor-Silberman \cite{Naor-Silberman} with very general results on random walks on groups acting on hyperbolic spaces of Maher-Tiozzo \cite{Maher-Tiozzo}.  In the i.i.d.\ labeling, the fact that edges are labeled independently at random has the consequence that, at arbitrarily large scales, the random walk on $G$ behaves like a random walk in one of the $\Omega_n$. By contrast, the graphical small cancellation condition ensures that the collection of words read on all the $\Omega_n$ is (at any given large enough scale) very sparse in the free group and, hence, in this model, the random walk on the group $G$ is very transversal to the random walks on the $\Omega_n$. In geometric terms, this is illustrated by the facts that, in the i.i.d.\ labeling, the $\Omega_n$ completely cover the Cayley graph of $G$ in the sense that coning off $\Cay(G,S)$ over the images of the $\Omega_n$ and their translates gives a bounded graph. On the other hand, in the graphical small cancellation labeling, coning off the Cayley graph of $G$ over the images of the $\Omega_n$ and their translates gives rise to a non-elementary hyperbolic graph \cite{GS-smallcanc}.

\medskip

In \cite{Naor-Silberman}, Naor-Silberman (building on an earlier work of Silberman \cite{Silberman}) prove that under strong spectral gap  properties for the sequence $\Omega_n$, the group $G$ almost surely satisfies fixed point properties for isometric actions on a large class  of metric spaces  (e.g.\ all CAT(0) spaces, all $L^p$-spaces for $p\in (1,\infty)$, \dots).  Although being very general, these metric spaces share a crucial \emph{local} property called uniform convexity. By contrast, our result concerns actions on spaces satisfying the \emph{large-scale} property of being Gromov hyperbolic. Our result thus naturally suggests the following question on a link between these two properties:
does there exist a group acting non-elementarily on a geodesic Gromov hyperbolic space such that for every $p\in (1,\infty)$, every isometric action on a $p$-convex geodesic metric space (or, more specifically, on an $L^p$-space) has a fixed point?
Indeed, answers to restrictions of our question to the subclasses of acylindrically hyperbolic groups or of graphical small cancellation groups would already be of great interest.

The following facts provide evidence towards a negative answer to our question: any Gromov hyperbolic group admits a proper action on some $L^p$-space \cite{Yu}. Moreover, all known examples of infinite groups for which all actions on $L^p$-spaces have fixed points are higher rank lattices \cite{BFGM}, or they come from Gromov's random construction. It is shown in \cite{Haettel} that higher rank lattices do not act non-elementarily on geodesic hyperbolic spaces, and our result shows that Gromov's groups coming from the i.i.d.\ labeling enjoy the same property.

\subsection*{Statement of the theorem}
Denote by $F_k$ the free group on the symmetric set of generators $S$ of size $2k$, and let $\Omega$ be a graph. A symmetric $F_k$--labeling $\alpha$ of $\Omega$ is a map from the edge set of $\Omega$ to $F_k$ so that $\alpha(e^{-1})=\alpha(e)^{-1}$ for every edge $e$. For $j\in\N$, we call $\mathcal A(\Omega,S^j)$ the set of symmetric $F_k$--labelings with values in $S^j$.

Consider $\Omega$ a disjoint union of finite connected graphs $\Omega_n$, i.e. $\Omega=\sqcup_{n\in\N}\Omega_n$, and endow $\mathcal A(\Omega,S^j)$ with the product distribution coming from the uniform distributions on the $\mathcal A (\Omega_n,S^j)$. For $\alpha\in \mathcal A (\Omega,S^j)$, define $G_\alpha$ to be the quotient of $F_k$ by all the words labeling closed paths in $\Omega$. For $n\in\N$, denote by $|\Omega_n|$ the cardinality of the vertex set of $\Omega_n$.

\begin{thm}\label{thm:main}
 Let $\Omega_n$ be a sequence of finite graphs, of vertex-degree between 3 and $d$ for some fixed $d\geq 3$. Assume $|\Omega_n|\to \infty$, and that there exists $C>0$ so that $\diam(\Omega_n)\leq C\girth(\Omega_n)$ for all $n$.
 
 Then for every $j\geq 1$ and almost every $\alpha\in \mathcal A(F_k,S^j)$, we have that $G_\alpha$ cannot act non-elementarily on any geodesic Gromov hyperbolic space.
\end{thm}

\begin{remark}
 Gromov's i.i.d.\ model also works for constructing monster groups as quotients of a given torsion-free non-elementary hyperbolic group. Observe that our result also holds in this setting because of the following argument. Let $H$ be a hyperbolic group generated by a set $S$, and for a labeling $\alpha$ of the $\Omega_n$ by $S^j$, define $H_\alpha$ as the quotient of $H$ by the relations as above. Then, in the notation of Theorem~\ref{thm:main}, $H_\alpha$ is a quotient of $G_\alpha$. Since the property of not admitting a non-elementary action on a geodesic hyperbolic space passes to quotients, this proves our claim.
\end{remark}

\subsection*{Outline of proof}

If a group $H$ acts non-elementarily on the hyperbolic space $X$, then by results of Maher-Tiozzo \cite{Maher-Tiozzo} random walks on $H$ make linear progress in $X$ with high probability. On the other hand, by results of Naor-Silberman \cite{Naor-Silberman}, a random walk on any one of the graphs $\Omega_n$ in the construction of a Gromov's monster $G_\alpha$ behaves like a random walk on $G_\alpha$ itself, up to a certain length. The hypothesis on the ratio between diameter and girth ensures that this length is at least some fixed fraction of the diameter of the graph. Hence, if $G_\alpha$ admitted a non-elementary action on a hyperbolic space $X$, by concatenating boundedly many trajectories of the random walk on one of the graphs $\Omega_n$ of that length (and using that generically the Gromov product at the concatenation point is small), one finds a path in $\Omega_n$ that maps to a path in $X$ whose diameter exceeds the diameter of (the image in $X$ of) $\Omega_n$. This is a contradiction.

\section{Notation and background}\label{sec:notation}

Let $\Lambda$ be a finite connected graph and denote by $F_k$ the free group on the symmetric set of generators $S$ of size $2k$. Let $(v_n)_{n\in\N}$ be the simple random walk on $\Lambda$ starting from the stationary measure.

We fix $j\geq1$ from now on. Let $\alpha_\Lambda\in \mathcal A(\Lambda,S^j)$, which we endow with the uniform distribution. Let $T$ be the Cayley tree of $F_k$. For a path $\gamma=(p_0,\dots,p_n)$ in $\Lambda$, there is a unique path\footnote{This is a slight abuse of language as two consecutive vertices $g_i$ and $g_{i+1}$ are joined by a path of length $\leq j$ and not by a single edge.} $(g_0,\dots,g_n)$ in $T$, starting at $1$, with the same edge-label as $\gamma$, and we denote by $\beta_{\alpha_\Lambda}(p_0,\dots,p_n)$ its endpoint. We denote by $\matP(\cdot)$ the probability of an event.

Define the distribution on $F_k$ $$\mu^n_{\Lambda,\alpha_\Lambda}(g)=\matP\bigl(\beta_{\alpha_\Lambda}(v_0,\dots,v_n)=g\bigr).$$

For $(w_n)$ the simple random walk on $F_k$ with respect to $S$ (i.e. the simple random walk on $T$ starting from 1), Naor-Silberman define $\overline\mu^n_{\Lambda,T}(g)$ as a certain convex combination $\sum_{l=0}^nP^n_\Lambda(l)\matP(w_{jl}=g)$, where $\sum_{l\leq n/6} P^n_\Lambda (l)\leq e^{-n/18}$.
By \cite[Proposition 7.4]{Naor-Silberman}, we have the following. Suppose that $\Lambda$ is a finite graph with vertex degrees between $3$ and $d<\infty$ and $\girth(\Lambda)\geq C\log |\Lambda|$, where $|\Lambda|$ denotes the cardinality of the vertex set. (Note that $\girth(\Lambda)\geq C\log |\Lambda|$ holds if the vertex degrees are bounded above by $d$ and $\diam(\Lambda)\leq \tilde C\girth(\Lambda)$ for $C$ depending only on $d$ and $\tilde C$.) Then there exist $C'>0$ depending only on the parameters $d,j,k,C$ (recall that $2k=|S|$) and for every $\lambda<1$ a function $\phi_\lambda:\N\to[0,1]$ going to $1$, depending on $\lambda,d,j,k$, such that:

$$\matP\bigl(\mu^n_{\Lambda,\alpha_\Lambda}(g)\geq \lambda\overline\mu^n_{\Lambda,T}(g) \text{ for all }n\leq C'\log|\Lambda|\bigr)\geq \phi_\lambda(|\Lambda|).\ \ \ \ (*)$$

\begin{remark}
 Naor-Silberman state a slightly weaker version of $(*)$, with $\lambda$ replaced by $1/2$, see \cite[Proposition 7.4]{Naor-Silberman}. However, their proof actually gives the statement above. In fact, replacing $1/2$ by $\lambda$ affects line --9 of the proof, where $\frac12 \epsilon(d,k,j)^q$ should be replaced by $(1-\lambda)\epsilon(d,k,j)^q$. Hence, in the formula at line --7, $\frac18$ should be replaced by $\frac{(1-\lambda)^2}{2}$. The condition on $C'$ in line --6 remains unchanged. Hence, the choice of $\lambda$ only affects the rate of convergence to $1$ of $\phi=\phi_\lambda$. 
\end{remark}

\medskip

\section{Proof of Theorem~\ref{thm:main}}

For the parameters $d,j,k,C$ appearing in Theorem~\ref{thm:main}, choose $C'$ as in $(*)$. By the Borel-Cantelli Lemma applied to $(*)$ (with $\Lambda=\Omega_m$ and $\alpha_\Lambda=\alpha|_{\Omega_m}$), 
for almost every $\alpha\in\mathcal A(\Omega,S^j)$ the following holds: for every $\lambda<1$ there exist infinitely many $m$ so that
$$\mu^n_{\Omega_m,\alpha|_{\Omega_m}}(g)\geq \lambda \overline\mu^n_{\Omega_m,T}(g) \text { for all } n\leq C'\log |\Omega_m|.\ \ \ \ (*)_m$$

Choose one such $\alpha$, and suppose by contradiction that $G_\alpha$ acts non-elementarily on the geodesic $\delta$--hyperbolic space $X$, so that $F_k$ does as well. Without loss of generality, we can assume that the orbit map $F_k\to X$ with respect to a fixed basepoint $x_0\in X$ is $1$--Lipschitz. For $g,h\in F_k$ denote $d_X(g,h)=d_X(gx_0,hx_0)$, and let Gromov products $(\cdot|\cdot)_{\cdot}$ be defined similarly, using $d_X$.

Recall that $(w_n)$ is the simple random walk on $F_k$ with respect to the generating set $S$. By \cite[Theorem 1.2, Lemma 5.9]{Maher-Tiozzo}, see also \cite{Mathieu-S}, we have some $\ell>0$ so that
\begin{equation*}
 \matP(d_X(1,w_n)\leq 6\ell n)\to 0,{\rm\ and}
 \end{equation*}
 \begin{equation*}
 \matP\left( (1|w_{2n})_{w_n}\geq \ell n/3\right)\to 0
\end{equation*}
as $n\to\infty$.

\begin{remark}
 Maher-Tiozzo state their results for actions on \emph{separable} geodesic hyperbolic spaces, but one can always pass to an action on a separable geodesic hyperbolic space. In fact, suppose that the countable group $H$ acts on the geodesic $\delta$--hyperbolic space $X$, for $\delta\geq 1$. Then one can take countably many $H$-orbits whose union is, say, $10\delta$--quasi-convex, by equivariantly choosing geodesics connecting pairs of points in a given orbit and considering a net in each such geodesic. By connecting $100\delta$--close pairs of points in the union of the given orbits, one obtains a graph $X'$ on which $H$ acts, and so that there exists an equivariant quasi-isometric embedding of $X'$ into $X$ (in particular, $X'$ is hyperbolic).
\end{remark}

\begin{remark}
 \cite[Lemma 5.9]{Maher-Tiozzo} is stated in a slightly different form than the one we use above, namely it gives control on the Gromov product $( (w_n^{-1}w_{2n})^{-1}| w_n)_1$. The proof only uses that $(w_n^{-1}w_{2n})$ is independent of $w_n$, so that it also applies to our case $(1|w_{2n})_{w_n}=(w_n^{-1}|w_{n}^{-1}w_{2n})_1$.
\end{remark}

Let $\lambda>0$ be very close to $1$, to be determined shortly. Suppose that $|\Omega_m|\geq e^{1/C'}$ is large enough that $\matP(d_X(1,w_n)\geq 6\ell n)\geq \lambda$ and $\matP\left((1|w_{2n})_{w_n}\leq \ell n/3\right)\geq \lambda$ for all $n\geq j\lfloor C'\log(|\Omega_m|) /6\rfloor$, where $C'$ is as in $(*)_m$. Suppose in addition that $(*)_m$ holds for the given $\lambda$.

By definition of $G_\alpha$, the graph $\Omega_m$ admits a j-Lipschitz map to $G_\alpha$, given as a label-preserving graph-homomorphism $\Omega_m\to\Cay(G_\alpha,S)$. Composition with the (1-Lipschitz) orbit map $G_\alpha\to G_\alpha\cdot x_0$ gives a j-Lipschitz map $\phi:\Omega_m\to X$. We let $\delta_X$ be the pull-back pseudo-metric on $\Omega_m$: $\delta_X(v,w)=d_X(\phi(v),\phi(w))$. Since $G_\alpha$ acts on $X$ by isometries, $\delta_X$ does not depend on the choice of label-preserving graph-homomorphism $\Omega_m\to\Cay(G_\alpha,S)$. Explicitly, if $(p_0,\dots,p_n)$ is a path in $\Omega_m$, then $\delta_X(p_0,p_n)=d_X(x_0,\beta_\alpha(p_0,\dots,p_n)\cdot x_0)$. 

In the sequel, Gromov products of vertices of $\Omega_m$ are defined using $\delta_X$.

As above, let $(v_n)$ be the simple random walk on $\Omega_m$ starting from the stationary distribution.
Fix $N=\lfloor C'\log(|\Omega_m|)\rfloor$ and an integer $\xi>8C/(\ell C')$. We claim that for all sufficiently large $m$ the following holds (setting $v_{-N}:=v_0$):
$$P:=\matP\bigl(\delta_X(v_{Ni},v_{N(i+1)})\geq \ell N \text{ and } (v_{N(i-1)}|v_{N(i+1)})_{v_{Ni}}\leq \ell N/3\ \forall \ i=0,\dots, \xi-1\bigr)>0.$$
In other words, for sufficiently large $m$ there exists a path $(v_0,\dots, v_{N\xi})$ in $\Omega_m$ satisfying the conditions in the formula above. If $\ell N$ is much larger than $\delta$, then by well-known facts about hyperbolic spaces, see e.g. \cite[Proposition 13]{Maher-S}, we have:
\begin{align*} 
\delta_X(v_0,v_{N\xi})&\geq \sum_{i=0}^{\xi-1} \Big(\delta_X(v_{Ni},v_{N(i+1)})-(v_{N(i-1)}|v_{N(i+1)})_{v_{Ni}} -(v_{Ni}|v_{N(i+2)})_{v_{N(i+1)}}-50\delta \Big) \\
&\geq \xi \ell N/4 >\diam(\Omega_m).
\end{align*}
This is a contradiction because the inequality above says that the image under the natural 1-Lipschitz map $\Omega_m\to X$ of the path we are considering has diameter larger than the diameter of $\Omega_m$.

We are left to prove the claim. Notice that 
\begin{align*}
\matP(\delta_X(v_{Ni},v_{N(i+1)})\geq \ell N)
&=\matP(\delta_X(v_{0},v_{N})\geq \ell N) \\ 
&=\mu_{\Omega_m,\alpha|_{\Omega_m}}(\delta_X(v_0,v_N)\geq \ell N) \\
&\geq \lambda\overline\mu^{N}_{\Omega_m,T}((\delta_X(v_{0},v_{N})\geq \ell N) \\
&=\lambda\sum_{l=0}^{N} P_{\Omega_m}^{N}(l)\matP(\delta_X(1,w_{jl})\geq \ell N) \\
&\geq\lambda\sum_{l=\lfloor N/6\rfloor+1}^{N} P_{\Omega_m}^{N}(l)\matP(\delta_X(1,w_{jl})\geq \ell N) \\
&\geq \lambda^2(1-e^{-N/18}),
\end{align*}
where in the first equality we used that the distribution of $(v_0,\dots,v_N)$ is the same as that of $(v_{Ni},\dots,v_{N(i+1)})$, which is true because the random walk starts from the stationary measure.

Similarly, $\matP((x_{N(i-1)}|x_{N(i+1)})_{x_{Ni}}\leq \ell N/3)\geq \lambda^2(1-e^{-N/18}).$
Hence:
$$1-P\leq 2 \xi(1- \lambda^2(1-e^{-N/18})),$$
which is less than $1$ if $\lambda>\sqrt{1-1/(2\xi)}$ and $N$ (which is a function of $|\Omega_m|$) is large enough.\qed

\medskip

\noindent{\bf Acknowledgment} We thank Goulnara Arzhantseva for helpful comments on the manuscript.

 \bibliographystyle{alpha}
 \bibliography{main}

\end{document}